\documentclass[openany,reqno,a4paper,12pt]{amsart}

\usepackage{amsmath,amssymb,amsthm}

\usepackage{mathrsfs}

\usepackage{mathtools}
\usepackage{commath}

\usepackage{thmtools}
\usepackage{thm-restate}

\usepackage{cases}
\usepackage{enumitem}
\setlist[enumerate,1]{label=(\roman*)}

\usepackage[pdftex, pdfborderstyle={/S/U/W 0}]{hyperref}
\hypersetup{
    colorlinks=true,
    linkcolor=magenta,
    citecolor=cyan,
}

\usepackage{cleveref}

\usepackage{etoolbox}

\usepackage{comment}

\usepackage[numbers]{natbib}

\numberwithin{equation}{section}

\ifdefined\thmcolor
\declaretheoremstyle[
  shaded={bgcolor=\thmcolor}
]{plain}
\else
\fi

\ifdefined\defcolor
\declaretheoremstyle[
  headfont=\normalfont\bfseries,
  bodyfont=\normalfont,
  shaded={bgcolor=\defcolor}
]{noital}
\else
\declaretheoremstyle[
  headfont=\normalfont\bfseries,
  bodyfont=\normalfont,
]{noital}
\fi

\declaretheorem[style=plain,numberwithin=section,name=Theorem]{theorem}

\declaretheorem[style=plain,sibling=theorem,name=Conjecture]{conjecture}

\declaretheorem[style=plain,numbered=no,name=Theorem]{theorem-n}
\declaretheorem[style=plain,numbered=no,name=Proposition]{proposition-n}
\declaretheorem[style=plain,numbered=no,name=Lemma]{lemma-n}
\declaretheorem[style=plain,numbered=no,name=Corollary]{corollary-n}
\declaretheorem[style=plain,numbered=no,name=Conjecture]{conjecture-n}
\declaretheorem[style=plain,numbered=no,name=Claim]{claim-n}
\declaretheorem[style=plain,numbered=no,name=Fact]{fact-n}
\declaretheorem[style=plain,numbered=no,name=Open Problem]{openproblem-n}
\declaretheorem[style=plain,numbered=no,name=Question]{question-n}

\declaretheorem[style=noital,numbered=no,name=Remark]{remark-n}
\declaretheorem[style=noital,numbered=no,name=Definition]{definition-n}
\declaretheorem[style=noital,numbered=no,name=Construction]{construction-n}
\declaretheorem[style=noital,numbered=no,name=Example]{example-n}

\undef{\set}
\DeclarePairedDelimiter{\set}{\lbrace}{\rbrace}

\undef{\emptyset}
\newcommand{\emptyset}{\varnothing}

\newcommand{\from}{\colon}

\undef{\abs}
\DeclarePairedDelimiterX{\abs}[1]
  {\lvert}{\rvert}{\ifblank{#1}{\,\cdot\,}{#1}}

\undef{\norm}
\DeclarePairedDelimiterX{\norm}[1]
  {\lVert}{\rVert}{\ifblank{#1}{\,\cdot\,}{#1}}

\DeclarePairedDelimiterX{\inner}[2]
  {\langle}{\rangle}{\ifblank{#1}{\,\cdot\,}{#1},\ifblank{#2}{\,\cdot\,}{#2}}

\DeclareMathDelimiter{\given}
  {\mathbin}{symbols}{"6A}{largesymbols}{"0C}

\DeclareMathOperator{\Prob}{\mathbb{P}}
\DeclarePairedDelimiterXPP{\prob}[1]
  {\Prob}{\lparen}{\rparen}{}
  {\renewcommand{\given}{\nonscript\;\delimsize\vert\nonscript\;\mathopen{}}#1}

\DeclareMathOperator{\Expec}{\mathbb{E}}
\DeclarePairedDelimiterXPP{\expec}[1]
  {\Expec}{\lparen}{\rparen}{}
  {\renewcommand{\given}{\nonscript\;\delimsize\vert\nonscript\;\mathopen{}}#1}

\DeclareMathOperator{\Var}{Var}
\DeclarePairedDelimiterXPP{\var}[1]
  {\Var}{\lparen}{\rparen}{}
  {\renewcommand{\given}{\nonscript\;\delimsize\vert\nonscript\;\mathopen{}}#1}

\DeclareMathOperator{\Cov}{Cov}
\DeclarePairedDelimiterXPP{\cov}[2]
  {\Cov}{\lparen}{\rparen}{}{#1,#2}

\newcommand{\RR}{\mathbb{R}}

\newcommand{\cC}{\mathcal{C}}
\newcommand{\cD}{\mathcal{D}}

\usepackage{geometry}
\usepackage{tikz}
\usetikzlibrary{quotes} 
\usepackage{standalone}

\newcommand{\corona}{$\cC_\text{odd}$}
\newcommand{\bin}{$\set{0,1}$}

\begin{document}

\title{On conjectures concerning the graph grabbing game}

\author{Lawrence Hollom}
\address{Department of Pure Mathematics and Mathematical Statistics (DPMMS), University of Cambridge, Wilberforce Road, Cambridge, CB3 0WA, United Kingdom}
\email{lh569@cam.ac.uk}



\begin{abstract}
  We consider two conjectures made in regard to the graph grabbing game, played on a vertex weighted graph.
  Seacrest and Seacrest conjectured in 2012 that the first player can win the graph grabbing game on any even-order bipartite graph.
  Eoh and Choi conjectured a strengthening of this in 2019, namely that the first player can win on any graph with no induced corona product of an odd cycle and a point.
  We provide a family of counterexamples to the latter conjecture, and propose a weaker conjecture in its place.
  We also show that the above two conjectures are equivalent when the vertex weights are all $0$ or $1$.
\end{abstract}

\maketitle


\section{Introduction}
\label{sec:intro}

The graph-grabbing game has two players, referred to as Alice and Bob.
They play on a graph $G$ with a weight function $w\from V(G)\to \RR_{\geq 0}$, alternately removing vertices of $G$ so that the remaining vertices induce a connected subgraph of $G$, with Alice going first.
Once all vertices are removed, the players compute their score, equal to the total weight of all vertices they removed.
Alice and Bob both aim to maximise their score, and we will say that Alice has a winning strategy, or simply that Alice \emph{wins} if she has a strategy guaranteeing that she collects at least half the total weight.

The game was introduced in Winkler's puzzle book \cite{winkler2003mathematical}, in which he asked for a proof that Alice could win when $G$ is an even-order path.
Micek and Walczak \cite{micek2011conjecture} then considered general trees of even order, proving that Alice can secure at least a quarter of the total weight.
This was in turn improved by Seacrest and Seacrest \cite{seacrest2012grabbing}, who supplied an elegant proof that Alice can always win on any tree of even order.
They furthermore made the following conjecture.

\begin{conjecture}[\cite{seacrest2012grabbing} Conjecture 5]
\label{conj:bipartite}
    Alice can win the graph-grabbing game played on a bipartite graph of even order.
\end{conjecture}

The conjecture remains open, but there has been partial progress; Egawa, Enomoto, and Matsumoto \cite{egawa2018graph} resolved the case of graphs they called $K_{m,n}$-trees (a copy of $K_{m,n}$ with trees attached to its vertices).
Eoh and Choi \cite{eoh2019graph} considered the situation where the weight function $w$ has codomain \bin. 
In this case, they were able to show that Alice wins whenever $G$ has no induced subgraph isomorphic to the corona product of a cycle and a point.
They also made the following conjecture.

\begin{conjecture}[\cite{eoh2019graph}]
\label{conj:corona}
    If \corona{} is the family of corona products of an odd cycle and a point, then Alice wins the graph grabbing game if and only if $G$ is \corona-free.
\end{conjecture}

More recently, Doki, Egawa, and Matsumoto \cite{doki2021graph} found, among other results, an example of a graph containing an induced \corona{} subgraph, yet on which Alice could still always win, ruling out one direction of \Cref{conj:corona}.

We give a family of graphs which are not \corona-free, yet there is some weight function for which Alice loses, resolving the other direction of \Cref{conj:corona} in the negative.
However, in \Cref{sec:binary-weights} we consider \Cref{conj:corona} in the special case of \bin-weighted graphs, and show that it is equivalent to \Cref{conj:bipartite} in this case.
Finally, we give a new conjecture on when Alice can guarantee to win at least half of the weight.


\section{Counterexamples}
\label{sec:counterexample}

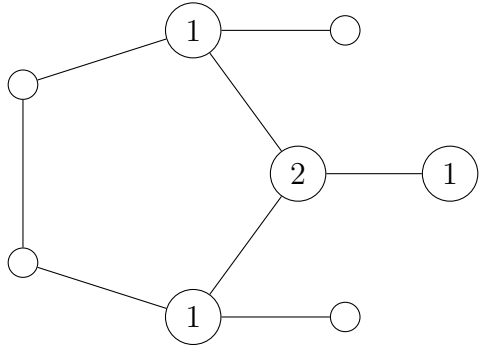
\begin{figure}[ht]
    \centering
    \begin{tikzpicture}[scale=2]
    \node[circle, draw=black] (A) at (1, 0) {2};
    \node[circle, draw=black] (B) at (0.31, 0.95) {1};
    \node[circle, draw=black] (C) at (-0.81, 0.59) {};
    \node[circle, draw=black] (D) at (-0.81, -0.59) {};
    \node[circle, draw=black] (E) at (0.31, -0.95) {1};

    \node[circle, draw=black] (El) at (1.31, -0.95) {};
    \node[circle, draw=black] (Al) at (2, 0) {1};
    \node[circle, draw=black] (Bl) at (1.31, 0.95) {};

    \draw[-] (A) -- (B);
    \draw[-] (B) -- (C);
    \draw[-] (C) -- (D);
    \draw[-] (D) -- (E);
    \draw[-] (E) -- (A);
    \draw[-] (A) -- (Al);
    \draw[-] (B) -- (Bl);
    \draw[-] (E) -- (El);
\end{tikzpicture}
    \caption{The simplest \corona-free connected graph on which Bob has a winning strategy. The numbers shown are vertex weights, and blank vertices have weight 0.}
    \label{fig:8-vtx-counterexample}
\end{figure}

The simplest \corona-free graph and weighting for which Bob can win is as shown in \Cref{fig:8-vtx-counterexample}.
It is not too difficult to see (though cumbersome to prove in writing) that Bob can indeed force a win on this weighted graph.

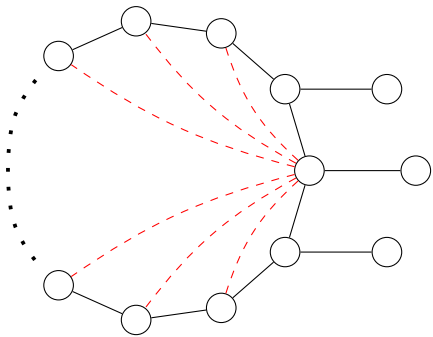
\begin{figure}[ht]
    \centering
    \begin{tikzpicture}[scale=2]
    \node[circle, draw=black] (A) at (1, 0) {};
    
    \node[circle, draw=black] (B) at (0.84, 0.54) {};
    \node[circle, draw=black] (C) at (0.42, 0.91) {};
    \node[circle, draw=black] (D) at (-0.14, 0.99) {};
    \node[circle, draw=black] (E) at (-0.65, 0.76) {};
    
    \node[circle, draw=black] (F) at (-0.65, -0.76) {};
    \node[circle, draw=black] (G) at (-0.14, -0.99) {};
    \node[circle, draw=black] (H) at (0.42, -0.91) {};
    \node[circle, draw=black] (I) at (0.84, -0.54) {};

    \node[circle, draw=black] (Il) at (1.51, -0.54) {};
    \node[circle, draw=black] (Al) at (1.7, 0) {};
    \node[circle, draw=black] (Bl) at (1.51, 0.54) {};

    \draw[-] (A) -- (B);
    \draw[-] (B) -- (C);
    \draw[-] (C) -- (D);
    \draw[-] (D) -- (E);

    \tikzstyle{dotted}= [dash pattern=on \pgflinewidth off 2mm]
    \draw[dotted, line width = 0.5mm] plot [smooth, tension=0.7] coordinates {(-0.8, 0.6) (-0.96, 0.28) (-0.96, -0.28) (-0.8, -0.6)};
    
    \draw[-] (F) -- (G);
    \draw[-] (G) -- (H);
    \draw[-] (H) -- (I);
    \draw[-] (I) -- (A);
    \draw[-] (A) -- (Al);
    \draw[-] (B) -- (Bl);
    \draw[-] (I) -- (Il);

    \draw[dashed, color=red] (A) to[bend left=15] (C);
    \draw[dashed, color=red] (A) to[bend left=12] (D);
    \draw[dashed, color=red] (A) to[bend left=10] (E);
    \draw[dashed, color=red] (A) to[bend right=10] (F);
    \draw[dashed, color=red] (A) to[bend right=12] (G);
    \draw[dashed, color=red] (A) to[bend right=15] (H);
\end{tikzpicture}
    \caption{The family of graphs $\cD$ which we conjecture must appear as an induced subgraph if Bob is able to win. The cycle has an odd number of vertices, and any subset of the dashed red edges may be present.}
    \label{fig:general-bob}
\end{figure}

This construction generalises to the family $\cD$ of graphs show in \Cref{fig:general-bob}.
Note that many members of $\cD$ (but not all members of any given order) contain smaller members of $\cD$ as induced subgraphs.
One may note that Bob can win on any graph in $\cD$ with the weighting generalised from that shown in \Cref{fig:8-vtx-counterexample} (so four vertices receive nonzero weight).
However, on larger graphs containing these as induced subgraphs, a more complex weighting may be necessary.

\section{\texorpdfstring{\bin}{\{0, 1\}}-weighted graphs}
\label{sec:binary-weights}

Eoh and Choi \cite{eoh2019graph} also discussed graphs where all vertex weights are in the set \bin, and conjectured that Alice could always win in \corona-free graphs in this case.
Here we prove the following result.

\begin{theorem}
    If Alice wins on every even-order bipartite graph with weights in \bin, then Alice wins on every even-order \corona-free graph with weights in \bin.
\end{theorem}

\begin{proof}
    Assume for contradiction that Bob wins on some even-order non-bipartite graph $G$ with weight function $w\from V(G)\to \set{0,1}$.
    We find a member of \corona{} induced in $G$.
    If $u$ is an optimal first play for Alice, and $v$ is an optimal play for Bob following $u$, then if $w(u)\geq w(v)$, then Bob must also win on $G-\set{u,v}$ with weights given by $w$.
    Let us assume that $G$ is a minimal graph (under removal of vertices) for which Bob wins.
    Thus for any non-cutvertex $u$ of $G$, we have $w(u)=0$ and there is some cutvertex $v$ of $G$ which is not a cutvertex of $G-u$ and has $w(v)=1$.
    So we see that every non-cutvertex of $G$ is a leaf.

    Now, $G$ is not bipartite, so take an induced odd cycle $x_1,\dotsc,x_r$ in $G$ (take an odd cycle and repeatedly pass to an odd sub-cycle if it is not induced).
    None of these vertices are leaves, so they must all be cutvertices of $G$.
    In particular, for each $x_i$ there is a vertex $y_i$ adjacent to $x_i$ such that there is no path from $y_i$ to $x_{i+1}$ in $G-x_i$ (subscripts taken modulo $r$).
    But then we see that $x_1,\dotsc,x_r,y_1,\dotsc,y_r$ induce an odd corona in $G$, with the points $x_i$ providing the cycle, and the points $y_i$ the leaves.
\end{proof}

Thus we see that in the \bin-weighted case, \Cref{conj:bipartite} and \Cref{conj:corona} are equivalent.


\section{Conclusion}
\label{sec:conclusion}

Having demonstrated that \Cref{conj:corona} does not hold in general, we propose the following weaker conjecture in its place.

\begin{conjecture}
\label{conj:new}
    If a graph $G$ of even order has no induced subgraph isomorphic to a member of the family $\cD$ (as shown in \Cref{fig:general-bob}), then Alice wins the graph grabbing game on $G$.
\end{conjecture}

We also note that \Cref{conj:bipartite} remains open, even for the case of \bin-weighted graphs.
Having attempted to generalise the approaches discussed in \Cref{sec:intro} which applied to simpler cases of this problem, it is the author's opinion that some significant new idea is needed to resolve these conjectures in their full generality.


\bibliographystyle{abbrvnat}  
\renewcommand{\bibname}{Bibliography}
\bibliography{main}


\end{document}